\documentclass[12pt]{amsart} 
\usepackage{amssymb}
\def\refeq#1{{\rm (\ref{eq #1})}}
\makeatletter
\let\over\@@over
\let\atop\@@atop
\let\above\@@above
\let\overwithdelims\@@overwithdelims
\let\atopwithdelims\@@atopwithdelims
\let\abovewithdelims\@@abovewithdelims
\makeatother

\let\cal\mathcal
\let\Bbb\mathbb

\def\calA{{\cal A}}
\def\calB{{\cal B}}
\def\calC{{\cal C}}
\def\calD{{\cal D}}
\def\calE{{\cal E}}
\def\calF{{\cal F}}
\def\calG{{\cal G}}
\def\calH{{\cal H}}
\def\calL{{\cal L}}
\def\calM{{\cal M}}
\def\calN{{\cal N}}
\def\calF{{\cal F}}
\def\calR{{\cal R}}
\def\calS{{\cal S}}
\def\calU{{\cal U}}
\def\calZ{{\cal Z}}

\def\<{\langle}
\def\>{\rangle}

\def\theequation{\arabic{section}.\arabic{equation}}
\renewcommand{\epsilon}{{\varepsilon}}

\newcommand{\beq}{\begin{eqnarray}}
\newcommand{\eeq}{\end{eqnarray}}
\newcommand{\R}{\mathbb{R}}
\newcommand{\N}{\mathbb{N}}
\newcommand{\Q}{\mathbb{Q}}
\newcommand{\C}{\mathbb{C}}

\newcommand{\nm}{{n,m}}
\newcommand{\nmk}{{n,m,k}}
\newcommand{\Gnmk}{{G_{n,m,k}}}

\newcommand{\FCbi}{{{\cal F}C_b^\infty}}
\newcommand{\FCbk}{{{\cal F}C_b^k}}
\newcommand{\Pdnm}{{P_{d_\nm}}}
\newcommand{\tozn}{{:z^n:\, }}
\newcommand{\toznl}{{:z^{n-1}:\, }}

\newtheorem{theorem}{\bf Theorem} 
\newtheorem{lemma}{\bf Lemma}
\newtheorem{prop}{\bf Proposition}
\newtheorem{corollary}{\bf Corollary}

\newtheorem{rem} {Remark}  

\newtheorem{de}{Definition}  
\renewcommand{\labelenumi} {\alph{enumi}}  
\renewcommand{\theequation}{\arabic{equation}}

\def\Ex{\Exaple\rm}
\newcommand{\ben}{\begin{equation}} \newcommand{\een}{\end{equation}}
\newcommand{\srs}{\stackrel{srs}{\longrightarrow}}
\newcommand{\pex}{{\PP}^x}  
\newcommand{\eex}{{\EE}^x}  
\newcommand{\indi}{\mbox{\rm 1\hspace{-2pt}\rule[0mm]{0.2mm}
{7.9pt}\hspace{ 1pt}}}
\begin{document}

\font\sy=cmsy9 scaled\magstep1
\font\cyr=wncyr10 scaled\magstep1
\font\cyrm=wncyr9 scaled\magstep0
\font\rc=eurm10 scaled\magstep1

\def\V{\Vert}
\def\v{\vert}

\def\o{\omega}
\def\O{\Omega}
\def\g{\gamma}
\def\G{\Gamma}
\def\var{\varphi}
\def\b{\beta}
\def\a{\alpha}
\def\n{\nabla}
\def\d{\delta}
\def\D{\Delta}
\def\Bscr {{\cal  B}}
\def\Ascr { {\cal  A}}
\def\Dscr { {\cal  D}}
\def\Fscr { {\cal  F}}
\def\Escr { {\cal  E}}
\def\Nscr { {\cal  N}}
\def\lscr { {\cal  l}}
\def\Mscr { {\cal  M}}
\def\Pscr { {\cal  P}}
\def\Sscr { {\cal  S}}
\def\Lscr { {\cal  L}}
\def\Hscr { {\cal  H}}
\def\Kscr { {\cal  K}}
\def\Tau {\frak T}

\def\hs {\hskip 10pt}
\def\leq {\leqslant}
\def\geq {\geqslant}
\def\A {{\frak A}_k}
\def\r {\restriction}
\def\rline{\hbox{{\rm I}\hskip -2pt {\rm R}}}  
\def\dbone{\hbox{1\hskip -3.5pt 1}}
\def\pr{^{\prime}}
\def\half{{1\over 2}}
\def\sig{\sigma}
\def\ol{\overline}
\def\eps{\varepsilon}
\def\lan{\langle}
\def\ran{\rangle}
\def\sig{\sigma}
\def\part{\partial}
\def\Lam{\Lambda}
\def\l{\lambda}
\def\rlinesm{\hbox{{\fiverm I}\hskip -1.5pt{\fiverm R}}}
\def\L{{\buildrel \circ \over L}}

\title[Strong uniqueness for certain Dirichlet operators]
{Strong uniqueness  for certain infinite dimensional Dirichlet operators
and applications to stochastic quantization}

\author{Vitali Liskevich}
\address{\hskip-\parindent
Vitali Liskevich\\
School of Mathematics\\
University of Bristol\\
Bristol BS8 1TW, UK}

\author{Michael R\"ockner}
\address{\hskip-\parindent
Michael R\"ockner\\
Fakult\"at f\"ur Mathematik\\
Universit\"at Bielefeld\\
 33501 Bielefeld, Germany} 
 \email{roeckner@mathematik.uni-bielefeld.de}

\begin{abstract}
Strong and Markov uniqueness problems in $L^2$ for Dirichlet operators on
rigged Hilbert spaces are studied. An analytic approach based on a--priori
estimates is used. The extension of the problem to the $L^p$-setting is
 discussed.
As a direct application essential self--adjointness and strong uniqueness in
$L^p$ is proved for the generator (with initial domain the bounded smooth
cylinder functions)  of the stochastic quantization process for Euclidean
quantum field theory in finite volume $\Lambda \subset \R^2$. 
 \end{abstract}

\subjclass{Primary: 47B25, 81S20; Secondary: 31C25, 60H15, 81Q10}
\keywords{Key words and phrases: Dirichlet operators, essential
self--adjointness, $C_0$--semigroups, generators, stochastic
quantization, Markov uniqueness, a--priori estimates}

\maketitle

\section{ Introduction}

The theory of Dirichlet forms is a rapidly developing field of modern analysis which has
intimate relationships with potential theory, probability theory, differential equations
and quantum physics. 
We refer to the  monographs \cite{BoHi}, \cite{D}, \cite{Fu} and \cite{MR92} where the 
theory of Dirichlet forms
with applications to different branches of analysis and probability theory is presented.
Though the abstract general theory is well developed, specific 
analytic questions remain open when one studies concrete situations. 
 In this paper we mainly discuss  Dirichlet forms and
corresponding Dirichlet operators on infinite-dimensional state spaces. In particular, 
we are concerned with the classical Dirichlet form of gradient type on a separable
Hilbert space with a probability measure. 
Initially, the form is defined on some ``minimal''
domain $\Dscr$.

The first analytic problem which arises when one
studies such forms is closability. This problem is well understood and necessary 
and sufficient conditions have been found (cf. \cite{AR1}). 
The operators associated
to the Dirichlet forms generate Markov semigroups (see e.g.
 \cite{D},\cite{Fu}, \cite{MR92}) on 
the corresponding $L^2$-space. If one assumes that the same domain $\Dscr$
 is contained in  the
domain of the generator, the next natural question which arises is 
whether the extension in $L^2$
of the generator restricted to $\Dscr$ is unique. There are at least two different statements
of the uniqueness problem in this context. 
The first is the so--called Markov uniqueness problem when one asks
whether the extension generating a Markov semigroup is unique. This problem is completely 
solved for the finite dimensional case in \cite{RZ1} where  Markov uniqueness was obtained
under the most general conditions. The situation is quite different in the 
infinite dimensional case.
We refer to \cite{ARZ93a}, \cite{ARZ93b} for the best results in this direction
 known so far.

One speaks of strong uniqueness for the Dirichlet form resp. operator
 when there is
only one lower bounded self-adjoint 
extension of the Dirichlet operator originally defined on $\Dscr$. 
As is well--known this is equivalent to essential selfadjointness. 
This problem was addressed
in many papers such as e.g.  
\cite{AKoR1},  \cite{AKRa}, \cite{AKoR2}, \cite{BKR96},
\cite{LSe1}, \cite{Shi}, 
\cite{W85} (see also the references therein). The results obtained
in these papers turn out to be incomparable, although they are all expressed in terms of 
conditions on the logarithmic derivative of the measure.

The aim of this paper is to investigate several ``uniqueness problems''
 connected with Dirichlet forms resp. Dirichlet 
operators on rigged Hilbert spaces. The main result (Theorem 1)
 concerns strong uniqueness in $L^2$.
The method we use is inherited from \cite{LSe1} (and also \cite{AKoR1},
 \cite{AKRa}, \cite{AKoR2}) and based upon an 
a--priori estimate of the cylindric
smooth solutions to the corresponding parabolic equations. This estimate is a generalization of
the estimate obtained in \cite{LSe1} to the case when the logarithmic derivative of the measure
contains two terms: one of them satisfies the conditions of \cite{LSe1} and the
 other is  modeled according to 
that in \cite{AKoR2} (see also \cite{AKoR1}, \cite{AKRa}).
The method of this paper is purely analytic, in
 contrast to that in \cite{AKoR2}.
It consists of the reduction of the problem
to estimates (independent of dimension) of gradients of solutions of the 
finite dimensional projections of the 
problem. We use these estimates to show that any  arbitrary semigroup generated by an 
extension of the original Dirichlet operator on $\Dscr$ 
can be approximated by the same approximation sequence, and therefore is unique. 
\medskip

While the main result of this paper is a (strict) generalization of \cite{LSe1},
it is, though similar in nature, still quite disjoint from those in
\cite{AKoR1}, \cite{AKRa}, \cite{AKoR2} w.r.t. applications. Nevertheless, 
our main
result can be applied to prove essential self--adjointness of the corresponding
Dirichlet operator in a situation where this, despite several attempts, had been 
an open problem for quite some time. That is  the so--called stochastic
quantization of field theory in finite volume. More precisely, here the
underlying space is a Sobolev space of distributions on an open bounded 
set in $\R^2$
and the measure is a two--dimensional Euclidean quantum field in  finite volume
with polynomial self--interaction. We refer to Section 5 below for the
(extensive list of the) corresponding literature and more details. We only mention
here that Markov uniqueness, however, 
had been shown already in \cite[Section 7]{RZ92}.
\medskip

Furthermore, it should be mentioned that 
 the  approach of the present paper to  the uniqueness problem naturally
extends
to the $L^p$-setting. The question is then whether the extension of the Dirichlet operator
on $\Dscr$ generating a $C_0$-semigroup on $L^p$ is unique.
For the first $L^p$-uniqueness result based upon the method of a--priori estimates we
refer to \cite{L1}. In this paper we obtain the uniqueness in $L^1$ (Theorem 2) which 
requires much simpler a--priori estimates than that needed for the proof of essential 
selfadjointness. As a consequence we derive a new approximation criterium for Markov
uniqueness (Corollary 1). We also discuss the uniqueness problem in $L^p$ for $p>1,~p\not=2$
generalizing the result from \cite{L1},
 although we do not present the details for the a--priori estimates needed
for this in the present 
paper, in order to avoid overloading the reader with additional technicalities.
But we emphasize that the corresponding
 results then also apply to the cases described in
Section 5 mentioned above. For further results on strong uniqueness in $L^p$ we
refer to \cite{Sta97}, where uniqueness results of ''perturbative type'' for $p=1$ are
proved, and to \cite{E97} where, in particular, strong uniqueness in $L^p$ in
the above mentioned situation of the stochastic quantization is also proved, but
only for $1\le p<2$. The latter case corresponds to our Theorems 
\ref{theorem 2} and \ref{theorem 4}.
\medskip

The organization of this paper is as follows. In Section 2 we present the framework and
the main uniqueness results. In Section 3 we derive the a--priori estimates for gradients
of the solutions to parabolic equations with smooth cylindric coefficients. In Section 4 
we give proofs of the uniqueness results. In Section 5 we discuss the said
applications.
Section 6 is devoted to the discussion of the
uniqueness problem in $L^p$.

\section{ Framework and main results}

Let $\Hscr_0$ be a separable real Hilbert space with the inner product
$(\cdot,\cdot)_0$ and norm $\v\cdot\v_0$. Let 
$$
\Hscr_+\subset\Hscr_0\subset\Hscr_-
$$ 
be 
a rigging of $\Hscr_0$ by the Hilbert spaces $\Hscr_+$ and $\Hscr_-$ 
with the assumption that 
the embeddings are dense, continuous and belong to the Hilbert--Schmidt class.
Without loss of generality we can suppose then that there exists a selfadjoint operator
$T=T^\ast\geq 1$ in $\Hscr_0$ with $\Dscr(T)=\Hscr_+$ such that $T^{-1}$ is 
 Hilbert--Schmidt. We refer to \cite[Chap. 6, \S 3]{BK} 
 for the details. We use the orthonormal  
system $(e_i)_{i=1}^\infty$ of eigenvectors of $T$ as a basis in $\Hscr_0$:  
$Te_i=\lambda_i e_i ~(i=1,2,3,\dots)$. For $e_i\in \Hscr_+$, $x\in \Hscr_-$, we
define $x_i:=\ _+(e_i,x)_-$, where $_+(\ ,\ )_-$ denotes the dualization between
$\calH_-$ and $\calH_+$. Clearly, the norms in $\Hscr_{\pm}$ can be
calculated as follows
$$
\v x \v_-^2=\sum_{i=1}^\infty \lambda_i^{-2} x^2_i, \qquad
\v x\v_+^2=\sum_{i=1}^\infty \lambda_i^{2}x^2_i.
$$
For $N\in\N$ define $P_N : \calH_- \to \calH_+\subset \calH_0$ by
$$
P_N x := \sum_{i=1}^N x_i\, e_i\ ,\ x\in \calH_-\ .
$$
Below we mostly identify the linear span of $\{e_1,\ldots,e_N\}$ with $\R^N$.

We will denote the space of $k$-times continuously
differentiable bounded mappings from $\Hscr_-$ into a Banach space $X$ by
$C_b^k(\Hscr_-,X)$ ($C_b^0(\Hscr_-,X)\equiv C_b(\Hscr_-,X)$).
$C^2_b(\Hscr_-,X)$  is a Banach space with the norm
$$
\V f \V_{C_b^2}=\sup_{x\in \Hscr_-}(\V f(x)\V_X+
\V f\pr(x)\V_{\Lscr(\Hscr_-,X)}+\V f^{\prime\prime}(x)\V_{
\Lscr(\Hscr_-,\Lscr(\Hscr_-,X))}),
$$
where $\Lscr(\Hscr_-,X)$ denotes the space of all bounded linear
operators from $\Hscr_-$ to $X$. 
When $X$ is the set of complex numbers $\Bbb C$ we identify $f\pr(x)\in
\Lscr(\Hscr_-,\Bbb C)$ with $\widetilde f\pr(x)\in \Hscr_+$ and
$f^{\prime\prime}(x)\in \Lscr(\Hscr_-,\Lscr(\Hscr_-,\Bbb C))$ with
the operator $\widetilde f^{\prime\prime}(x)\in  
\Lscr(\Hscr_-,\Hscr_+)$. In this case 
$$
_+(\widetilde f\pr(x),\var )_-=f\pr(x)\var,~~~\ \ 
_+(\widetilde f^{\prime\prime}(x)\var,\psi)_-=(f^{\prime\prime}(x)\var)
\psi,~~\var, \psi\in\Hscr_-.
$$
So we make the convention that $\nabla f:=f^\prime=\widetilde f^\prime$ and 
$f^{\prime\prime}=\widetilde f^{\prime\prime}$ for $f\in C_b^2(\Hscr_-)$. 

By $\FCbk$, $k\in\N \cup\{\infty\}$, 
let us denote  the set of all functions $f$ on $\Hscr_-$
such that there exist $N\in \Bbb N$, 
$\{\phi_1,\dots,\phi_N\}\subset \calH_+$
and $G \in C_b^k(\R^N)$ such that
$$
f(x)=G (\ _+(\phi_1,x)_-, \dots,\ _+(\phi_N,x)_-), \ x\in \Hscr_-.
$$

Let $\nu$ be a probability measure defined on the $\sigma$-algebra
$\Fscr$ of Borel subsets of $\Hscr_-$ with $\hbox{supp}\, \nu=\Hscr_-$. 
Assume that $\nu$ has a 
logarithmic derivative in the sense  that there exists an 
$\Fscr$-measurable mapping
$\b:\Hscr_-\mapsto\Hscr_-$ and $\forall q\in \Hscr_+,\ 
 \b_q(\cdot):=\ _+(q,\b(\cdot))_- \in L^2(\Hscr_-,\nu)$ 
and the following integration by parts formula holds
\begin{equation}\label{eq 2.1}
\int_{\Hscr_-}\n_if(x)d\nu(x)=-\int_{\Hscr_-}\b_i(x)f(x)d\nu(x)\ \mbox{for all
  $f\in \FCbi$}.
\end{equation}
where $\b_i:=\ _+(\b,e_i)_-$. We also use the following notations
$$
\n_if:=(\n f,e_i)_0,\qquad \D f:=Tr_{\Hscr_0}f^{\prime\prime}.
$$
$\V\cdot\V_p$ is the norm in $L^p(\Hscr_-,\nu)\equiv L^p,~~
\lan \cdot,\cdot\ran$ is the inner product in $L^2$, and $\<\ \ \>$ denotes
expectation w.r.t. $\nu$.

We introduce the Dirichlet form in $L^2$ as a closure of the form
$$
\Escr(f,g)\equiv\int_{\Hscr_-}(f\pr(x),g\pr(x))_0\ d\nu = 
\lan (\n f, \n g)_0\ran,~~f,g\in \FCbi .
$$
Under the stated conditions the form is closable (see \cite{AR1}).
The integration by parts formula implies that
$$\lan \n_i f, g\ran=\lan f,-(\n_i+\b_i)g\ran~~\forall f,g\in C_b^1.$$

 Let $A$ be the operator associated with the form $\Escr$.
Then $A=A^\ast \geq 0$. The semigroup
$(e^{-tA},t \geq 0)$ is positivity preserving (i.e.,
$e^{-tA}[L^2_+] \subset L^2_+:=\{f \in L^2: f \geq 0~~\hbox{a.e.}\}$)  and 
       $L^\infty$-contractive
(i.e., $\V e^{-tA} f \V_\infty \leq \V f \V_\infty,~~\forall f \in L^2 \cap
L^\infty$). We also have that $e^{-tA}\,1=1$ for all $t>0$, since $A1=0$.
Recall that such a semigroup is called Markov semigroup and its generator is called
Markov generator. 
It is a standard fact that this semigroup defines a family of $C_0$-semigroups of contractions
in $L^p,~~p\in [1,\infty)$:
$$e^{-tA_p}:=(e^{-tA}\restriction [L^2\cap L^p])^\sim_{L^p\rightarrow L^p}.$$

In what follows we use a stronger condition on the logarithmic derivative,
namely, we assume that $\b_i\in L^2(\Hscr_-)$ which enables us to identify the action
of the operator $A$ :
$$
A\restriction C_b^2=-(\D+\, _-(\b,\n\cdot)_+) \restriction C_b^2.
$$

The operator $A$ associated with the Dirichlet form is called the
{\em Dirichlet operator}.

In the following definitions we define different types of the uniqueness problem for the
operator $A$.

\begin{de}\label{definition 1}
A linear set $\Dscr\subset\Dscr(A)$ is called a domain
of strong uniqueness of $A$ if there exists only one selfadjoint extension
of $A\restriction\Dscr$ in $L^2$.
\end{de}
We note that according to a theorem by J. von Neumann this 
is equivalent to the essential selfadjointness of $A\restriction\Dscr$.

\begin{de} \label{definition 2}
A linear set $\Dscr\subset\Dscr(A)$ is called a domain
of Markov uniqueness of $A$ if there exists only one extension
of $A\restriction\Dscr$  which is a Markov generator in $L^2$.
\end{de}
It is clear that if $\Dscr$ is a domain of strong uniqueness of $A$ it is also a domain
of Markov uniqueness, but not vice versa. (There are counterexamples even if
$\calH_- = \calH_0 =\calH_+ =\R^1$, cf. \cite{E97}.) 

We also extend the strong uniqueness problem to the $L^p$-setting.
\begin{de}\label{definition 3}
Let $p\geq 1$. A linear set $\Dscr\subset\Dscr(A_p)$ is called a domain
of $L^p$-strong uniqueness of $A$ if all extensions $B$
of $A\restriction\Dscr$ in $L^p$ such that $-B$ is a generator of a $C_0$-semigroup
in $L^p$ coincide with $A_p$.
\end{de}

According to a result by W.Arendt 
(cf. \cite{N} Theorem A-II, 1.33, p.46) 
the latter
is equivalent to the fact that $\Dscr$ is a core of the operator $A_p$.

The following simple observation gives a link between the $L^1$-strong uniqueness and
the Markov uniqueness.
\begin{prop}\label{proposition 1}
Let $A$ be a Markov generator. Let $\Dscr\subset\Dscr(A)\cap\Dscr(A_1)$.
 Then if $\Dscr$ is a domain of $L^1$-strong uniqueness for $A$,
it is also a domain of Markov uniqueness.
\end{prop}
{\it Proof}
Let $f\in\Dscr$. Then by the definition of $A_1$ we have $e^{-tA_1}f=e^{-tA}f$ and
$${1\over t}(f-e^{-tA_1}f)={1\over t}(f-e^{-tA}f)={1\over t}\int_0^te^{-tA}Afds
={1\over t}\int_0^te^{-tA_1}Afds.$$
Passing to the limit with respect to $t\downarrow 0$ we get $A_1f=Af$.
Now the result follows easily.\qed
\bigskip

Now we are ready to formulate the main results of the paper. The first result concerns
the strong uniqueness problem in $L^2$. In \cite{LSe1} it was proved that if
$\b:\Hscr_-\mapsto\Hscr_0$ and $\v\b\v_0\in L^4(\Hscr_-,d\nu)$ then 
$-\D -\, _-(\beta,\nabla\cdot)_+ \restriction C_b^2$ is essentially selfadjoint.
However, in the infinite dimensional case these conditions become quite
restrictive in applications, in particular,
 because of the condition $\b(x)\in \Hscr_0$ for $\nu$--a.e.  $x\in \Hscr_-$. 
The authors of  \cite{AKoR2} (using related  analytic ideas,
 but  also stochastic techniques in an essential way) 
derived another sufficient condition for essential selfadjointness,
assuming the existence of the derivative of $\b$ and that it satisfies 
certain one-sided estimates.
The following result is a strict generalization of \cite{LSe1} invoking, in
addition, also a one sided estimate (but involving the $|\ |_+$--norm rather 
than the $|\ |_-$--norm as in \cite{AKoR2}).
 The primary idea is to consider $\b$ as a sum of  mappings
$\alpha:\Hscr_-\mapsto\Hscr_0$ and $\delta:\Hscr_-\mapsto\Hscr_-$. 
The former satisfies
the conditions of \cite{LSe1} whereas the derivative of the
latter will satisfy a one-sided estimate. Technically this is achieved by 
means of an  a--priori estimate (see Section 3). Though this result remains 
quite disjoint from those in  \cite{AKoR1}, \cite{AKRa}, \cite{AKoR2} 
we are nevertheless
able to apply it to new examples not covered by those papers. The
most striking one is to show essential selfadjointness of the Dirichlet operator
corresponding to the stochastic quantization of finite volume quantum
fields. This is presented in Section 5 below, to which we refer for details and references.
\begin{theorem}\label{theorem 1}
Let $\beta = \alpha + \delta$, $|\alpha |_0 \in L^4(\calH_-,\nu)$, $|\delta |_-
\in L^2 (\calH_-,\nu)$. Suppose that there exist $\delta^m : \calH_- \to
\calH_-$, $m\in \N$, such that
\begin{itemize}
\item[(i)] $\delta^m = \delta^m\circ P_{N_m}$ for some $N_m\in\N$ and
  $\delta_j^m \restriction \R^{N_m}\in C_b^1 (\R^{N_m})$ with globally 
H\"older
  continuous first order derivatives for all $1\le j\le N_m$ (where $\delta^m_j
  := \, _+(e_j, \delta^m)_-$).
\item[(ii)] $|\delta - \delta^m|_- \to 0$  in $L^2 (\calH_-,\nu)$ as
  $m\to \infty$.
\item[(iii)]  There exists $c_+ \in \R$ such that for all $m\in\N$, 
$$
(\Lambda_{\delta^m}(x)y,y)_+\le c_+\ |y|^2_+ \ \forall x,y\in\R^{N_m},
$$
where $\Lambda_{\delta^m}(x):= (\nabla_i\delta_j^m (x))_{1\le i,j\le N_m}$.
\item[(iv)]
If $\alpha \neq 0$, suppose, in addition, that $\delta_j=\delta_j\circ P_j$ with
$\delta_j\restriction \R^j \in C^1(\R^j)$ for all $j\in\N$ and that  there
exists $\epsilon_0\in (0,1)$, $c(\epsilon_0)\in\R$ such that for all $N\in\N$
\begin{eqnarray*}
&&\lefteqn{
\left\<\left( (\nabla_i\delta_j\circ P_N)_{1\le i,j\le N}\, w,w\right)_0
\right\> } \\
  &&
  \le (1-\epsilon_0)\sum_{i,j =1}^N \<\nabla_{\!j}\, w_i\circ P_N,\nabla_{\!j}\, w_i\circ P_N\> + c(\epsilon_0)
  \<(w\circ P_N,w\circ P_N)_0\>
\end{eqnarray*} 
\vspace{-2mm}

\noindent
for all $w_1,\ldots,w_N\in C_b^1(\R^N)$ and $w:=  (w_1,\ldots,w_N)$. \\
\end{itemize}
Then
$(-\Delta -\, { }_-\!(\beta,\nabla\cdot)_+ \restriction \FCbi)$ is
essentially self-adjoint in $L^2(\calH_-,\nu)$.

\end{theorem}

\begin{rem}  Consider the following condition:
\begin{itemize}
\item[(iii')]
 $(\Lam_{\d^m}(x)y,y)_-\leq c_-\v y\v^2_- \ \ \forall x,y \in\R^{N_m}$.
\end{itemize}
\noindent
To compare (iii) and (iii') one can rewrite them in the form
$$(iii)\ \ \ (T\Lam_{\d^m}(x)T^{-1}z,z)_0\leq c_+\v z\v_0^2 \ \ 
\forall x,z\in \R^{N_m}. 
$$
and
$$(iii')\ \ \ (T^{-1}\Lam_{\d^m}(x)Tz,z)_0\leq c_-\v z\v_0^2 \ \ 
\forall x,z\in \R^{R_N}.
$$

Therefore (iii) and (iii') coincide if $\Lambda_{\d^m}=\Lambda^\ast_{\d^m}$.
\end{rem}

The next theorem gives a criterium for the $L^1$-strong uniqueness of the
 operator $A$.

\begin{theorem}\label{theorem 2}
Let $\b=\a+\d,~~\v \a\v_0, \v \d\v_-\in L^2(\Hscr_-,\nu)$.
Suppose that there exists a sequence of mappings $(\d^m)_{m\in\Bbb N},\ 
\d^m:\Hscr_-\mapsto\Hscr_-, m\in \Bbb N$ such that

(i) $(\d^m)_{m\in\Bbb N}$ satisfies condition (i) in Theorem 
\ref{theorem 1}.

(ii) $\v\d-\d^m\v_-\longrightarrow 0$ in $L^1(\Hscr_-,d\nu)$ as
$m\rightarrow\infty$.

(iii) There exists a constant $c_+\in \Bbb R$ such that for all $m\in\N$
$$
(\Lam_{\d^m}(x)y,y)_+\leq c_+\v y\v^2_+ \ \ \forall x,y\in \R^N.
$$

\noindent Then the operator $(\D + \, _-(\b,\n\cdot)_+\restriction
\FCbi)$ has  a unique extension which generates a $C_0$-semigroup
on $L^1(\calH_-,\nu)$.
\end{theorem}

\begin{corollary}\label{corollary 2.3}
Let the conditions of Theorem \ref{theorem 2}
 be satisfied. Then $\calF C_b^\infty$ is a domain of 
Markov uniqueness for the operator $A$.
\end{corollary}

The proofs of Theorems 1 and 2 will be given in Section 4 after we derive 
a--priori estimates in Section 3. These estimates are the core of the method.

\section{A--priori estimates}

The aim of this section is to obtain a--priori estimates for solutions of 
parabolic equations on $\R^d$. So here 
$\Hscr_0, \Hscr_+,\Hscr_-$ will be just
$\ \R^d$ endowed with the inner products
 $(y,z)_+=\sum_{j=1}^d\mu_j^2 y_j z_j,~(y,z)_-=\sum_{j=1}^d\mu_j^{-2}
 y_j z_j$ and the usual Euclidean product for $\Hscr_0$. The measure $\nu$ above
 is correspondingly now a measure on $\R^d$.

Let \(u(t,x)=u\) be the solution of the Cauchy problem
\begin{equation}\label{eq 3.1}
 \left\{
 \begin{array}{cc}
     {\partial  \over \partial t} u =\D u +(b,\n u)_0\\
          u(\cdot,0)=f(\cdot),
 \end{array}
 \right.
\end{equation}
where $f\in C_0^\infty( \R^d),\ b\in C^1_b( \R^d, \R^d)$ with globally H\"older
continuous first order derivatives.
The operator $\Ascr:=\D+(b,\n\cdot)_0$ generates a $C_0$-semigroup on
$C_\infty(\Bbb R^d)$ which can be extended to  a $C_0$-semigroup on $L^p(\Bbb
R^d)$ ($:=L^p(\R^d,dx)$, where $dx$ denotes Lebesgue measure) and for 
 the solution of \refeq{3.1} we have $u(t,\cdot)=e^{t\calA}f(\cdot)$, $t>0$,
 and  $\V u\V_\infty\leq \V f\V_\infty$. Furthermore, we have  $u(t,\cdot)\in
C_b^2 (\R^d)\cap C^3(\R^d)$ (even with globally H\"older continuous second
order derivatives and locally H\"older continuous third order derivatives; cf
\cite[Theorems 9.2.3 and 8.12.1]{Kr96}).
\bigskip

Let us introduce the derivative of the mapping $b:\Bbb R^d\mapsto \Bbb R^d$
as the linear operator $\Lambda_b$ whose matrix is $(\Lambda_b)_{ij}=(\n_i b_j)$.

\begin{prop}\label{proposition 2}
Let $u$ be the solution to \refeq{3.1} with $f\in C_0^\infty(\R^d)$.
Suppose that $b$ is as above satisfying, in addition, that there exists $c_+\in
\R^+$ such that 
$$
(\Lambda_b(x) y,y)_+\leq c_+\ (y,y)_+\ \ \forall x,y\in \R^d.
$$
Then
\begin{equation}\label{eq 3.2}
\V \v \n u(t,x)\v_+\V_\infty \leq e^{c_+ t}\V \v \n f
\v_+\V_\infty.
\end{equation}
\end{prop}
\medskip

\noindent{\bf Proof.}
Let $w_i:=\n_iu$ be the derivative of $u$ in direction $e_i$. Denote the inner product in $L^2(\Bbb
R^d)$ by $\lan \cdot,\cdot \ran$ and the integral with respect to Lebesgue
measure $dx$ by $\lan\cdot\ran$
({\it only in this proof}).
Differentiating equation \refeq{3.1} in the direction $e_i$ we get
$$
{d\over dt}w_i-\D w_i-(b, \n w_i)_0-(\Lambda_b w)_i=0,
$$
where $(\Lambda_b w)_i=\sum_{j=1}^d(\n_ib_j)w_j.$ 
Multiplying both sides of the last equality
by $\mu_i^2w_i\v w \v_+^{p-2}$, $p>4$, 
after integration with respect to $dx$ and summation with
respect to $i$ from $1$ to $d$ we have

\begin{eqnarray*}
{1\over p}{d\over dt}\V \v w\v_+\V_p^p
+\sum_{i,j}\lan \n_j w_i,\v w\v_+^{p-2}\mu_i^2\n_j w_i\ran+
\sum_{i,j}\lan w_i\mu_i^2\n_j w_i,\n_j\v w\v_+^{p-2}\ran\\ 
 -\sum_{i,j}\lan(\n_i b_j)w_j,\mu_i^2 w_i\v w\v_+^{p-2}\ran-
\sum_i\lan(b,\n w_i)_0,\mu_i^2 w_i\v w\v_+^{p-2}\ran=0.
\end{eqnarray*}
The last term is equal to
\begin{eqnarray*}
&&-\half \<(b,\nabla |w|_+^2)_0\ |w|_+^{2(p/2-1)}\> = -{1\over p}\<(b,\nabla
|w|_+^p)_0\> =  {1\over p} \<|w|_+^p \ \mbox{div}\, b\>\ .
\end{eqnarray*}
Therefore, we obtain the equality
\begin{eqnarray}\label{eq 3.3}
{1\over p}{d\over dt}\V \v w\v_
+\V_p^p+4{p-2\over p^2}\V | \n (\v w
\v_+^{p/2})|_0 \V_2^2+\sum_{j=1}^d\lan (\n w_j,\n w_j)_+,\v w
\v_+^{p-2}\ran
\nonumber\\
=\lan (\Lambda_b\ y,y)_+\ran-{1\over p}\lan |w|_+^p \mbox{div}\ b\ran,
\end{eqnarray}
where $y:=w\v w \v_+^{p/2-1},~~\V\cdot \V_p$ is the norm in $L^p(\R^d)$.
From \refeq{3.3} and the assumption  it follows that
$$
{1\over p}{d\over dt}\V \v w\v_+\V_p^p\leq c_+ \V \v w\v_+\V_p^p+
{k\over p}\V \v w\v_+\V_p^p,~~\hbox{~where~}~k:=\V div\
b\V_\infty
$$
or
$$\V \v w(t)\v_+\V_p\leq \V \v w(0)\v_+\V_p\ e^{(c_+ +{k\over
p})t}.$$
Since $\v w \v_+(0)=\v \n f \v_+$, passing to the limit $p \rightarrow \infty$ we get
$\V \v w\v_+\V_\infty \leq \V \v \n f \v_+\V_\infty\ e^{ct}$.\qed
\bigskip

From now on as in Section 2 the norms $\|\ \|_p$  again denote the $L^p$--norms
w.r.t $\nu$.

\begin{prop}\label{proposition 3}
Let $\b=\alpha+\d,\ b=\alpha^1+ \delta^1$ 
with $\alpha, \delta, \alpha^1, \delta^1 :
\R^d\to \R^d$. Suppose that
$\v\alpha\v_0\in L^2(\R^d,\nu),\ \v\d\v_-\in L^1(\R^d,\nu)$ and that the 
condition of Proposition \ref{proposition 2} is
satisfied. Let $u$ be the solution of \refeq{3.1}. Then

\begin{eqnarray*}
&&\|u\|_2^2+ 
\int_0^t\V\v\n u\v_0\V_2^2\ ds\\
& \leq & t\V f\V_\infty^2\V\v\alpha-\alpha^1\v_0
\V_2^2+2c_+^{-1}(e^{c_+t} -1)\|\, | \nabla f|_+\, \|_\infty\, 
\V f\V_\infty\, 
\V\v\d-\delta^1 \v_-\V_1+\V f\V_2^2.
\end{eqnarray*}
\end{prop}
\medskip

\noindent{\bf Proof.}
Multiplying \refeq{3.1} by $u$ and integrating w.r.t. $\nu$ we obtain 
after integration by parts that
$$
{1\over 2}{d\over dt}\V u \V_2^2+\V\v\n u\v_0\V_2^2=\lan(b-\b,\n u)_0, u\ran.
$$
Estimating the r.h.s. as follows
$$\v\lan(b-\b,\n u)_0,u\ran\v\leq {1\over 2}\V\v\n u\v_0\V_2^2+ {1\over 2}
\V f\V_\infty^2\V\v\alpha-\alpha^1\v_0\V_2^2+
\V f\V_\infty\V\v\d-\delta^1\v_-\|_1\, \V\v\n u\v_+\V_\infty$$
after integration with respect to $t$ and using Proposition 
\ref{proposition 2} 
one completes the proof.\qed

\begin{prop}\label{proposition 4}
Let $u$ be the solution to \refeq{3.1}.
In addition to the conditions of Proposition \ref{proposition 3}
 assume that  $\v \a \v_0\in L^4(\R^d,\nu)$, 
$|\delta|_- \in L^2(\R^d,\nu)$, 
$\delta \in C^1(\R^d,\R^d)$, and
 that there exist $\epsilon_0\in (0,1)$, $c(\epsilon_0)\in\R_+$ such that
\begin{equation}\label{eq 3.4}
\lan (\Lam_\d w,w)_0\ran\leq
(1-\eps_0)\sum_{j=1}^d\lan(\n_jw,\n_jw)_0\ran+c(\eps_0)\lan
(w,w)_0\ran
\end{equation}
for all  $w=(w_j)\in C_b^1(\R^d,\R^d)$. 
Then there exists $C(\eps_0)\in\R_+$  (depending  only on $\eps_0$) such 
that  
\begin{eqnarray}\label{eq 3.5}
&& \int_0^t \V \v \n u(s) \v_0\V_4^4ds \leq C (\eps_0) 
\{t\V f\V_\infty^4(\V\v\alpha\v_0\V_4^4+\V\v \alpha^1 \v_0\V_4^4
+ \|\, |\alpha -\alpha^1|_0\|_2^2) \nonumber \\
&& +{1\over 2c_+} (e^{2c_+t}-1)\V f\V_\infty^2\V\v \n f\v_+\V_\infty^2
\V\v\d-\delta^1\v_-\V_2^2\nonumber\\
&& +2c_+^{-1}(e^{c_+t}-1) \V f \V_\infty^3\V\v \n f\v_
+\V_\infty
\V\v\d-\delta^1\v_-\V_1
+\V f\V_\infty^2\V\v\n f\v_0\V_2^2
+c\V f \V_2^2\ \|f\|_\infty^2\}.\nonumber\\
&&
\end{eqnarray}
\end{prop}

Before proving Proposition \ref{proposition 4}
 we prove several lemmas (all of them under the
conditions of Proposition \ref{proposition 4}).
\begin{lemma}\label{lemma 1}  Let $w:=\nabla u$. Then
\begin{equation}\label{eq 3.6}
\V {du\over dt}\V_2^2+{d\over dt}\V\v w\v_0\V_2^2\leq\V(\b-b,
w)_0\V_2^2.
\end{equation}
\end{lemma}
\medskip

\noindent {\bf Proof.}
Using the equation and integrating by parts we have
$$
\lan (b, w)_0,{du\over dt}\ran=\lan {du\over dt}-  div\, w,{du\over dt}\ran=
\V {du\over dt}\V_2^2+\lan (w,{dw\over dt})_0\ran+\lan(\b, w)_0,{du\over dt}
\ran .
$$
Rewriting this in the form
$$
\V {du\over dt}\V_2^2+{1\over 2}{d\over dt}\V\v w\v_0\V_2^2=\lan (b-\b,
w)_0,{du\over dt}\ran \le \half \|(\beta -b,w)_0\|_2^2 +\half \|{du\over
  dt}\|_2^2, 
$$
we get the result.\qed

\begin{lemma}\label{lemma 2}
Let $u$ be the solution to (2), $w=\n u$. Then
\begin{eqnarray}\label{eq 3.7}
\nonumber
&&{1\over 2}{d\over dt}\V\v w\v_0\V_2^2+\sum_{i,j}\lan \n_i w_j,\n_i w_j\ran
= - \sum_i\lan (\alpha,\n w_i)_0,w_i\ran
+\lan (\Lambda_{\delta}w,w)_0\ran\\
&&+\lan (\delta-b,w)_0,{du\over dt}\ran
+\lan (\delta-\delta^1, w)_0,(\alpha-2\alpha^1, w)_0\ran\\
\nonumber
&&-\lan (\alpha^1, w)_0,
(\alpha, w)_0\ran + \|(\alpha^1,w)_0\|_2^2 + \V (\delta-\delta^1, w)_0 \V_2^2.
\end{eqnarray}
\end{lemma}
\medskip

\noindent{\bf Proof.}
Differentiating the equation in direction $e_i$ we get
$$
{dw_i\over dt}=\D w_i+\n_i(\alpha^1, w)_0+(\delta^1, \n w_i)_0+(\n_i\delta^1)
\cdot w.
$$
Multiplying scalarly by $w_i$ in $L^2(\R^d,\nu)$ and summing over $i$ we 
have the equality
$$
\displaylines{
{1\over 2}{d\over dt}\V\v w\v_0\V_2^2+\sum_{i,j}\lan \n_i w_j,\n_i
w_j\ran+\sum_{i,j} 
\lan(\alpha_j+\delta_j)\n_j w_i,w_i\ran+\lan (\alpha^1,
w)_0,\sum_i(\n_i+\b_i)w_i\ran\cr
-\sum_i\lan(\delta^1,\n w_i)_0,w_i\ran-\lan(\Lambda_{\delta^1}w,w)_0\ran=0.\cr
}  
$$
Again using the equation we rewrite the last equality in the form
$$\displaylines{
{1\over 2}{d\over dt}\V\v w\v_0\V_2^2+\sum_{i,j}\lan \n_i w_j,\n_i w_j\ran+
\sum_i\lan (\alpha, \n w_i)_0,w_i\ran+
\lan (\alpha^1, w)_0,(\alpha-\alpha^1, w)_0\ran\cr
+\lan (\alpha^1, w)_0,{du\over dt}\ran
+\lan (\alpha^1,
w)_0,(\delta-\delta^1, w)_0\ran-
\lan (\Lambda_{\delta^1}w,w)_0\ran+\sum_i\lan(\delta-\delta^1, \n
w_i)_0,w_i\ran=0.\cr
}$$
To finish the proof of the lemma one should observe that
$$\displaylines{
\sum_i\lan(\delta-\delta^1, \n
w_i)_0,w_i\ran=\lan (\Lambda_{\delta^1}w,w)_0\ran-\lan (\Lambda_{\delta}w,w)_0\ran\cr
-\lan(\delta-\delta^1, w)_0,{du\over dt}\ran-\lan(\delta-\delta^1,
w)_0,(\alpha-\alpha^1,w)_0\ran-\V(\delta-\delta^1, w)_0\V_2^2.\qed\cr}$$
\begin{lemma}\label{lemma 3}
Let $u$ be the solution to (2), $w:=\n u$. Then
$$
\V \v w\v_0\V_4^4\leq 16 \V f\V_\infty^2(\half \V(b-\b, w)_0\V_2^2-{1\over 4}
{d\over dt}\V\v w\v_0\V_2^2+ \sum_{i,j}\V \n_iw_j\V_2^2).
$$
\end{lemma}
\medskip

\noindent{\bf Proof.}
Using  equation \refeq{3.1} and integrating by parts  we obtain 
$$\displaylines{\V\v w \v_0\V_4^4=\sum_{i,j}\lan w_i,w_iw_j^2\ran=
-\sum_{i,j}\lan u, w_j^2\n_i w_i\ran-2 \sum_{i,j}\lan u,w_iw_j\n_i
w_j\ran-\sum_{i,j}\lan u,\b_iw_iw_j^2\ran\cr
=-\lan u, \v w \v_0^2{du\over dt}\ran-2 \sum_{i,j}\lan u,w_iw_j\n_i
w_j\ran+\lan u, \v w \v_0^2(b-\b, w)_0 \ran\cr
\leq \V u\V_\infty\V\v w\v_0\V_4^2\V{du\over dt}\V_2+\V u\V_\infty\V\v
w\v_0\V_4^2\V(b-\b, w)_0\V_2+
2\V u\V_\infty\sum_{i,j}\V w_iw_j\V_2\V\n_i w_j\V_2\cr
}$$
$$\displaylines{
\leq {1\over 2}\V\v w\v_0\V_4^4+ \V f\V_\infty^2\V{du\over dt}\V_2^2+
\V f\V_\infty^2\V(b-\b, w)_0\V_2^2\cr
+2\V f\V_\infty(\sum_{i,j}\V w_iw_j\V_2^2)^{1/2}(\sum_{i,j}\V \n_i w_j\V_2^2)^{1/2}\cr
\leq {3\over 4}\V\v w\v_0\V_4^4+\V f\V_\infty^2(\V(b-\b, w)_0\V_2^2-{d\over
dt}\V\v w\v_0\V_2^2)\cr
+\V f\V_\infty^2\V(b-\b, w)_0\V_2^2+4\V f\V_\infty^2\sum_{i,j}\V \n_i
w_j\V_2^2,
\cr}
$$
where we used that $ab\le \frac{1}{4} a^2 + b^2$ and Lemma \ref{lemma 1} in
the last step.  Now the assertion follows.
\qed
\bigskip

\noindent{\bf Proof of Proposition \ref{proposition 4}.}
Let us first estimate the first term of the right hand side of \refeq{3.7}:
\begin{eqnarray*}
\v\sum_i\lan (\alpha,\n w_i)_0,w_i\ran\v &\leq &\lan \v\alpha\v_0\v
w\v_0(\sum_{i,j}\v\n_iw_j\v^2)^{1/2}\ran\\
&\leq & {\eps_0\over 2}\sum_{i,j}\V \n_iw_j\V_2^2+{1\over 2\eps_0}\V
\v\alpha\v_0\cdot\v  w\v_0\V_2^2.
\end{eqnarray*}
The term in \refeq{3.7} containing ${du\over dt}$ is estimated by Lemma 
\ref{lemma 1}  as follows
\begin{eqnarray*}
&&\v\lan ( \delta -b , w)_0,{du\over dt}\ran\v\leq {1\over 2}\V {du\over
dt}\V_2^2+\half \V(\delta - b, w)_0 \V_2^2\\
& \leq& {1\over 2}\V(b-\b, w)_0\V_2^2-{1\over 2}{d\over dt}\V\v w\v_0\V_2^2+
\V(\delta -\delta^1, w)_0\V_2^2 + \| (\alpha^1, w)_0\|_2^2\ .
\end{eqnarray*}
Therefore, using that $ab \le \half a^2 +\half b^2$ for 
estimating the rest of the right hand side of \refeq{3.7}  we get by assumption
\refeq{3.4} for any $\eps,\eps_1>0$ from \refeq{3.7} that
\begin{eqnarray*}
&&{d\over dt}\V\v w\v_0\V_2^2+{\eps_0\over 2} (\sum_{i,j}\V \n_iw_j\V_2^2
+\half \| (b-\beta ,w)_0\|_2^2)\\
&\leq & {1\over 2\eps_0}\V\v\alpha\v_0\v w\v_0\V_2^2+
\tilde C(\eps_0)(\V(\alpha^1, w)_0\V_2^2+\V(\alpha, w)_0\V_2^2+\V(\d-\delta^1,
w)_0\V_2^2) +c(\eps_0)\V\v w\v_0\V_2^2\\
&\leq & {\eps\over 2\eps_0}\V\v w\v_0\V_4^4+{1\over
8\eps_0\eps}\V\v\alpha\v_0\V_4^4+\eps_1 \V\v w\v_0\V_4^4+ {\tilde C(\eps_0)
\over 2\eps_1}(\V\v \alpha^1\v_0\V_4^4+\V\v \alpha\v_0\V_4^4)\\
&& +\tilde C(\eps_0)\V(\d-\delta^1, w)_0\V_2^2+
c(\eps_0)\V\v w\v_0\V_2^2.
\end{eqnarray*}
where $\tilde C(\eps_0) \in \R_+$ only depends on $\eps_0$. Estimating the 
left hand side of  
the last inequality from below by  Lemma \ref{lemma 3}  one obtains 
the desired result after integration with respect to $t$ applying Propositions
\ref{proposition 2} and \ref{proposition 3}, and  
properly choosing $\eps$ and $\eps_1$
(for instance, $\eps:={\eps_0^2\over 64\V f\V_\infty^2}$,
$\eps_1={\eps_0\over 128 \V f\V_\infty^2}$).
\qed

\section{Proofs of the uniqueness results}

In this section we prove the results formulated in Section 2. Our strategy 
for the uniqueness in $L^2$ and in $L^1$
is the same.
Namely, we take an arbitrary extension
of the original operator defined on $\Fscr C_b^\infty$, which generates a $C_0$-semigroup in
$L^2$ ($L^1$ respectively). Then we construct the approximation sequence which converges
to the generated semigroup in $L^2$ ($L^1$ respectively). This implies the uniqueness of the
semigroup, and therefore the uniqueness of the extension.
The same approach is extendable to all $L^p$-spaces (see Section 6 for
more details).
\bigskip

\noindent{\bf Proof of Theorem 1.}  Let $f\in \FCbi$ and let  $N\in\N$, $G\in
C_b^\infty  (\R^N)$ such that
$$
f(x)=G(e_1(x),\ldots,e_N(x)),\quad x\in\calH_-\ .
$$
For every $n\in \N$ there exist
$$
\alpha_1^n,\ldots ,\alpha_{a_n}^n\in C_b^\infty (\R^{a_n})
$$
such that  for $\alpha^n \in C_b^\infty (\calH_-,\calH_-)$ defined by
$$
\alpha^n(x):=\sum_{i=1}^{a_n} \alpha_i^n(e_1(x),\ldots,e_{a_n}(x))\ e_i,\quad
x\in \calH_-\ ,
$$
we have 
$$
|\alpha -\alpha^n|_0 \longrightarrow 0 \ 
\mbox{in $L^4(\calH_-,\nu)$ as $n\to \infty$}.
$$
We may assume that $a_n\ge N$ for all $n\in\N$. For $m,n\in\N$ define
$$
d_{n,m}:=\max (m,a_n,N_m)
$$
where $N_m$ is as in assumption (i) of the theorem. We now apply the results of
Section 3 with $b$ replaced by $b^{n,m}:= P_{d_\nm}\circ (\alpha^n
+\delta^m)\restriction \R^{d_\nm}$ and $\nu$ replaced by $\nu\circ
P_{d_\nm}^{-1}$. For $k\in\N$ let $\chi_\nmk\in C_0^\infty (\R^{d_\nm})$, 
such that for all $n,m\in\N$ the maps $\chi_\nmk$, $\nabla \chi_\nmk$, 
$k\in\N$, are
uniformly bounded and such that $\chi_\nmk (x)=1$ provided $|x|\le k$. Define
$\Gnmk : = \chi_\nmk\cdot G$, where $G$ is considered as a function on
$\R^{d_{n,m}}$.
Then the solution of the Cauchy problem on $\R^{d_\nm}$
\begin{eqnarray*}
&&{\partial u_\nmk \over \partial t} =  \Delta u_\nmk + (b^\nm , \nabla 
u_\nmk)_0\\
&& u_\nmk (\ ,0)=\Gnmk
\end{eqnarray*}
is given by a $C_0$--semigroup on $C_\infty (\R^{d_\nm})$, i.e.,
$$
u_\nmk (t,\cdot) = e^{-t\calA_\nm} \Gnmk \ ,\quad t>0,
$$
where
$$
\calA_\nm = -\Delta - (b^\nm , \nabla\,\cdot)_0\ \mbox{on $C_0^\infty
  (\R^{d_\nm})$.} 
$$
Let now $B$ with domain $D(B)$ be an arbitrary lower bounded self--adjoint
extension of $(-\Delta -\ {}_-(\beta,\nabla\cdot)_+)\restriction \FCbi$ 
on $L^2(\calH_-,\nu)$. 
It is an easy exercise to see that $\calF C_b^2 \subset D(B)$
and that
$$
B = -\Delta - {}_-\!(\beta,\nabla\cdot)_+\ \mbox{on $\calF C_b^2$},
$$
in particular, $(e^{-t\calA_\nm}\Gnmk)\circ P_{d_\nm}\in D(B)$, $t\ge 0$. Since
$t\mapsto e^{-t\calA_\nm}\Gnmk$ is continuously differentiable from $\R_+$ to
$C_\infty (\R^{d_\nm})$, so is $t\mapsto (e^{-t\calA_\nm}\Gnmk)\circ P_{d_\nm}$ from
$\R_+$ to $L^2(\calH_- ,\nu)$. Therefore, 
\begin{eqnarray*}
&& -{d\over ds} \bigg (e^{-(t-s)B}\big(e^{-s\calA_\nm}\Gnmk\big)\circ
P_{d_\nm}\bigg)\\
&=& e^{-(t-s)B} \bigg(\big(\calA_\nm\ e^{-s\calA_\nm}\Gnmk \big)\circ P_{d_\nm}
\bigg)\\   
&& - B\ e^{-(t-s)B} \bigg(\big(e^{-s \calA_\nm}\Gnmk \big)\circ P_{d_\nm}\bigg)\\
&=& e^{-(t-s)B}\big(-\Delta - (b^\nm \circ P_{d_\nm},\nabla \cdot)_0 -B\big)
\left(e^{-s\calA_\nm}\Gnmk\right)\circ P_{d_\nm}\ .
\end{eqnarray*}
Hence we have justified that the classical Duhamel formula applies in our case,
since the  above implies that for all $t>0$ and $f_k:=\Gnmk \circ P_{d_\nm}$
\begin{eqnarray*}
&&e^{-tB} \ f_k - (e^{-t\calA_\nm}\Gnmk)\circ P_{d_\nm}\\
&=&\int_0^t e^{-(t-s)B}\ {}_-\!(\beta -b^\nm\circ P_{d_\nm}, \nabla
(e^{-s\calA_\nm} \Gnmk)\circ P_{d_\nm})_+\ ds\ .
\end{eqnarray*}
Hence, if $\gamma\in\R$ such that $B\ge \gamma$, then 
\begin{eqnarray*}
&&\sup_{k\in\N} \|e^{-tB} f_k - (e^{-t\calA_\nm}\Gnmk)\circ P_{d_\nm}\|_2\\
&\le & e^{t\gamma}\sup_{k\in\N} \int_0^t \|\, {}_-\!\left((\alpha^n+\delta^m)\circ P_{d_\nm}
-\beta, \nabla (e^{-s\calA_\nm}\Gnmk)\circ P_{d_\nm}\right)_+\|_2\ ds\\
&\le & e^{t\gamma}
\bigg[ \| \, |\alpha^n -\alpha|_0\, \|_4\quad  \sup_{k\in\N} \int_0^t  
\|\,|\nabla u_\nmk (s,P_{d_\nm} (\cdot))|_0\|_4\ ds\\
&&\qquad +  \|\, |\delta^m -\delta\|_- \|_2\quad  \sup_{k\in\N}\int_0^t 
\|\,|\nabla u_\nmk
(s,P_{d_\nm})|_+\,\|_\infty\ ds\bigg].
\end{eqnarray*}
Here we used the fact that, since $d_\nm \ge \max (a_n, N_m)$,
 both $\alpha^n\circ
P_{d_\nm}=\alpha^n$ and $\delta^m\circ P_{d_\nm} =\delta^m$. We want to show that
for all $t>0$
\begin{equation}\label{eq 4.8}
\overline{\lim_{n\to\infty}}\ \overline{\lim_{m\to \infty}}
\sup_{k\in\N} \|e^{-tB}f_k - (e^{-t\calA_\nm} \Gnmk)\circ \Pdnm
\|_2 =0. 
\end{equation}
Applying Proposition \ref{proposition 2} we see that as $m\to\infty$
the second term converges
to zero by assumption (ii) for each fixed $n$. We note that the assumption in
Proposition \ref{proposition 2} is indeed satisfied with $b=b^\nm$ for fixed
$n$ and a constant $c_+$ independent of $m$, since $\alpha^n\restriction
\R^{d_\nm}$ has bounded continuous derivatives and because of assumption 
(iii). If $\alpha =0$, \refeq{4.8} follows. If $\alpha\neq 0$, we apply
Proposition \ref{proposition 4} to see that the first term converges to zero,
too, if $n\to \infty$. We note that by assumption (iv) 
 we can indeed apply Proposition \ref{proposition 4}, since the
logarithmic derivative $\beta^\nm$ of $\nu\circ P_{d_\nm}^{-1}$ (as
 is easily checked) satisfies the equation
$$
\beta_j^\nm\circ \Pdnm = E_\nu [\beta_j\mid \sigma (\Pdnm)] = E_\nu
[\alpha_j\mid \sigma (\Pdnm)] +\delta_j
$$
for all $j\in \{1,\ldots,d_\nm\}$. Here 
$E_\nu [\, \cdot\, \mid \sigma (\Pdnm)]$
denotes the conditional expectation of $\nu$ given the $\sigma$--algebra 
$\sigma (\Pdnm)$ generated by  $\Pdnm$. Of course, the constant $c_+$ in the
right hand side of inequality \refeq{3.5} depends on $n$ and $m$. But since we
first take $m\to \infty$ it disappears so when afterwards taking $n\to\infty$,
the right hand side of \refeq{3.5} stays bounded.

Since $e^{-tB}f_k\longrightarrow e^{-tB}f$ in $L^2(\calH_-,\nu)$ as
$k\to\infty$, equality
\refeq{4.8} implies that $e^{-tB}f$, $t>0$, is independent of which extension
$B$ we took in the first place. Since $\FCbi$ is dense in
$L^2(\calH_- ,\nu)$, it follows that $(e^{-tB})_{t\ge 0}$ and hence that $B$ is
uniquely determined. Thus, the theorem is proved.\qed

\bigskip

\noindent{\bf Proof of Theorem 2.}
The idea of the proof is similar to that of Theorem 1.
We use the same approximating semigroups $e^{-t\Ascr_{m,n}}$ as in the proof of Theorem 1.

Let $B$ be an arbitrary  extension of 
$(-\D-\, _-(\b,\n\cdot)_+\restriction \calF C_b^\infty)$ such that $B$ is 
the generator of 
a $C_0$-semigroup in $L^1(\Hscr_-,\nu)$.

Now let $f$, $f_k$, $\Gnmk$ be as in the proof of Theorem \ref{theorem 1}.  
Then again $(e^{-s\calA_\nm}G_\nmk)\circ P_{d_\nm} \in D(B)$ and 
using Duhamel's formula we obtain the estimate
\begin{eqnarray*}
&&\V e^{-tB}f_k-(e^{-t\Ascr_{n,m}}\Gnmk\circ P_{d_\nm})\V_1\\
&\leq& e^{t\g} \{
\int_0^t\V (b^{n,m}\circ P_{d_\nm}-\b, \n 
(e^{-s\Ascr_{n,m}}\Gnmk)\circ P_{d_\nm})_0\V_1ds\}\\
&\leq & e^{t\g}\{
\V \v \a^n-\a\v_0\V_2 \sup_{k\in\N} \int_0^t \V\v \n u_{n,m,k}(s)\v_0\V_2ds
+\V \v \d-\d^m\v_-\V_1 \sup_{k\in\N}\int_0^t\V\v \n u_{n,m,k}(s)\v_+\V_\infty ds\}.
\end{eqnarray*}
Now we proceed as above, but using Proposition \ref{proposition 3} instead of
Proposition \ref{proposition 4} to show that

$$
\overline{\lim_n}\overline{\lim_m}\sup_{k}
\V e^{-tB}f_k-(e^{-t\Ascr_{n,m}}\Gnmk)\circ P_{d_\nm}\V_1=0.
$$
Therefore, the extension $B$ is unique which implies the assertion of
the theorem.\qed 

\section{Application to the stochastic quantization of field theory in finite
  volume} 

In this section we present our main application. After the programme of
stochastic quantization of field theories was initiated by Parisi and Wu in
\cite{PaWu81} and after its implementation for Euclidean quantum fields with
polynomial interaction in finite volume by Jona--Lasinio and Mitter in the
pioneering work \cite{J-LMi85},  there has been an enormous number of follow-up
papers on this subject (see e.g.  
\cite{AKR96}, \cite{AR89a}, \cite{AR89b},
 \cite{AR91}, \cite{ARZ}, \cite{ARZ93a}, \cite{ARZ93b}, 
\cite{BoCMi87}, \cite{DaPTu96}, \cite{GaGo95},\cite{HuKa96},
\cite{J-LMi90}, \cite{J-LSe91}, \cite{MikRoz97}, \cite{Mi86}, \cite{RZ92}, 
\cite{RZ96}. 
But the question
whether the corresponding Dirichlet operator restricted to $\FCbi$ is
essentially self--adjoint remained open. Below, we shall show that we can 
settle
this problem positively as a simple application of Theorem \ref{theorem 1}
above. We note that Markov uniqueness in this case was already shown in
\cite[Section 7]{RZ92}. We use the same notation as in the latter paper, but
for completeness we now recall the complete framework. 

Let $\Lambda$ be an open rectangle in $\R^2$. Let $(-\Delta +1)_N$ be the
generator of the following quadratic form on $L^2(\Lambda, dx): (u,v)\mapsto
\int_\Lambda \<\nabla u, \nabla v\>_{\R^2} dx + \int_\Lambda uv\ dx$ with
$u,v\in \{g\in L^2(\Lambda ,dx)\mid \nabla g \in L^2 (\Lambda, dx)\}$ (where
$\nabla$ is in the sense of distribution).  Let $\{e_n\mid n\in\N\}\subset
C^\infty (\bar \Lambda)$ be the (orthonormal) eigenbasis of $(-\Delta +1)_N$ and
$\{\lambda_n\mid n\in\N\}$ ($\subset ]0,\infty[$) the corresponding eigenvalues
(cf. \cite[p.~226]{RS78}), i.e., we consider Neumann boundary conditions. Define
for $\alpha \in \R$
\begin{equation}\label{eq 5.1}
H_\alpha :=\{u\in L^2(\Lambda, dx)\mid \sum_{i=1}^\infty
\lambda_n^\alpha\<u,e_n\>^2_{L^2(\Lambda,dx)} <\infty\}
\end{equation}
equipped with the inner product
\begin{equation}\label{eq 5.2}
\<u,v\>_{H_\alpha}:= \sum_{n=1}^\infty \lambda_n^\alpha
\<u,e_n\>_{L^2(\Lambda,dx)}\ \<v,e_n\>_{L^2(\Lambda,dx)}. 
\end{equation}
Clearly, we have that
\begin{equation}\label{eq 5.3}
H_\alpha \left\{\begin{array}{l}
\mbox{completion of $C^\infty (\bar\Lambda)$ w.r.t $\|\ \|_{H_\alpha}$ if
  $\alpha \ge 0$}\\
\mbox{completion of $C^\infty_0 (\Lambda)$ w.r.t $\|\ \|_{H_\alpha}$ if
  $\alpha \le 0$}
\end{array}
\right.
\end{equation}
(cf. \cite[p.~79]{LM72} for the latter).

Fix $\delta >0$. Since $\sum_{i=1}^\infty \lambda_n^{-1-\delta}<\infty$, we have
applying \cite[Theorem 3.2]{Y89} (i.e., the Gross--Minlos--Sazonov theorem) 
with
$H:=L^2(\Lambda,dx)$, $\|\cdot\|:=\|\cdot\|_{H_{-\delta}}$, $A_1:=(-\Delta
  +1)_N^{-\delta/2}$ and $A_2:=(-\Delta +1)^{-1/2}_N$, 
that there exists a unique mean zero Gaussian 
probability measure 
$\mu$ on  $\calH_- := H_{-\delta}$ (called {\em free field} on
$\Lambda$; see \cite{N73}) such that
\begin{equation}\label{eq 5.4}
\int_E l(z)^2\ \mu(dz) = \|l\|_{H_{-1}}^2\quad \mbox{for all $l\in
  \calH_-'=H_\delta$.} 
\end{equation}
Clearly, supp $\mu=E$. 

\begin{rem}\label{remark 5.1}
In \refeq{5.4} we have realized the dual of $H_{-\delta}$ as $H_\delta$ using
as usual the chain
\begin{equation}\label{eq 5.5}
H_\alpha\subset H_0=L^2(\Lambda,dx)\subset H_{-\alpha},\quad \alpha \ge 0.
\end{equation}
\end{rem}

Let $h\in L^2(\Lambda ,dx)$, $n\in\N$, and define $:z^n:\, (h)$ as follows 
(cf., e.g. \cite[Sect.~8.5]{GJ86}): fix $n\in\N$ and let $H_n(t)$, $t\in\R$, be
the $n$-th Hermite polynomial, i.e., $H_n(t)=\sum_{m=0}^{[n/2]} 
(-1)^m\alpha_{nm} t^{n-2m}$, with $\alpha_{mn}=n! / [(n-2m)!\, 2^m m!]$. Let
$d\in C_0^\infty(\R^2)$, $d\ge 0$, $\int d(x)\ dx=1$, and $d(x)=d(-x)$ for each
$x\in\R^2$. Define for $k\in\N$, $d_{k,x}(y):=2^{2k}d(2^k(x-y))$; 
$x,y\in\R^2$. Let $z_k(x):={}_{+}\!(d_{k,x},z)_-$, $z\in \calH_-$,
$x\in\Lambda$, and set
\begin{equation}\label{eq 5.6}
:z_k^n:\,(x)\, := c_k(x)^{n/2} H_n(c_k(x)^{-1/2}z_k(x)),
\end{equation}
 where $c_k(x):=\int z_k(x)^2\ \mu(dz)$. Then it is known that \\
$:z_k^n:\, (h)\, :=\int :z_k^n:(x)\, h(x)\ dx
\mathop{\longrightarrow}\limits_{k\to  \infty}
 :z^n:\, (h)$ both in every $L^p(\calH_-,\mu)$, $p\ge 1$, and for
$\mu$--a.e. $z\in \calH_-$ (cf., e.g., \cite[Sect.~3]{R86} for the latter). The
function $z\mapsto \lim \sup_{k\to\infty} :z_k^n:(h)$  
is then a $\mu$--version of $:z^n:\, (h)$. From now on
$:z^n:\, (h)$ shall denote this particular version.

Now fix $N\in\N$, $a_n\in\R$, $0\le n\le 2N$ with $a_{2N}<0$ and define
$$
V(z):= \sum_{n=0}^{2N} a_n\, \tozn (1_\Lambda),\quad z\in E,
$$
 where $1_\Lambda$ denotes the indicator function of $\Lambda$. Let
\begin{equation}\label{eq 5.7}
\varphi:=\exp(-\half V).
\end{equation}
Then $\varphi  >0$  $\mu$--a.e. and $\varphi\in L^p(\calH_-,\mu)$ for all 
$p\in
[1,\infty[$ (cf. e.g. \cite[Sect.~5.2]{Si74} or \cite[Sect.~8.6]{GJ86}). Set
$$
\nu:=\varphi^2\cdot\mu.
$$
We want to apply Theorem \ref{theorem 1} to $\nu$, $\calH_- := H_{-\delta}$ 
(as above), $\calH_0:= H_\alpha$, where
\begin{equation}\label{eq 5.8}
\alpha > \max (0,1-\delta/2)
\end{equation}
and $\calH_+ := H_{\delta+2\alpha}$. Clearly, then since by \refeq{5.8}, in
particular, $\alpha > 1-\delta$, the embeddings
$$
\calH_+ \subset \calH_0\subset \calH_-
$$
are Hilbert--Schmidt. We note that thus $\calH_+ = (-\Delta+1)^\alpha
(\calH_-')$ and that $\{\lambda_j^{-\alpha/2}e_j\mid j\in\N\}$,
$\{\lambda_j^{\delta/2}e_j\mid j\in\N\}$ are orthonormal bases of $\calH_0$,
$\calH_-$ respectively. Define $\beta:=\alpha+\delta : \calH_-\to \calH_-$,
where $\delta :\calH_-\to\calH_-$ and $\alpha : \calH_-\to\calH_0$ are for
($\nu$--a.e.) $z\in \calH_-$ defined by
\begin{eqnarray*}
&&\delta(z):=-\sum_{j=1}^\infty \lambda_j^{-\alpha -{\delta \over 2} +1}\,
{}_{H_{-\delta}}\!(z,e_j)_{H_\delta}\ \lambda_j^{\delta/2}e_j\\
&&\alpha(z):= -\sum_{j=1}^\infty \lambda_j^{-\alpha/2}\ \sum_{n=1}^{2N} n\ a_n\
\toznl (e_j)\ \lambda_j^{-\alpha/2} e_j\ .
\end{eqnarray*}
We have that $|\delta|_-$ , $|\alpha|_0^2\in L^2(\calH_-,\nu)$. Indeed,
\begin{eqnarray*}
\int |\delta|_-^2\ d\nu &=& \int \sum_{j=1}^\infty \lambda_j^{-2\a -\d +2}\
{}_{H_{-\delta}}(z,e_j)_{H_\delta}^2\ \varphi^2(z)\ \mu(dz)\\
&\mathop{\leq}\limits_{\refeq{5.4}} & \sum_{j=1}^\infty
\lambda_j^{-2\a -\d +2}\, 3^{1/4}\
\l_j^{-1} \left( \int \varphi^4\ d\mu\right)^{1/2}\\
&<& \infty,  
\end{eqnarray*}
since $2\a + \d -1<\infty$ by \refeq{5.8}, and
\begin{eqnarray*}
\left( \int |\alpha|_0^4\ d\nu\right)^{1/4} &=& \left(\int\left[\sum_{j=1}^\infty
  \l_j^{-\a} \left(\sum_{n=1}^{2N} n\ a_n\ \toznl (e_j)\right)^2\right]^2\
  \varphi^2(z) \ \mu(dz)\right)^{1/4}\\
&\le & \sum_{n=1}^{2N} n\ a_n\ \left(\int\left[\sum_{j=1}^\infty
  \l_j^{-\a}\left(\toznl (e_j)\right)^2\right]^2\ \varphi^2(z)\
  \mu(dz)\right)^{1/4}\\
&\le & \sum_{n=1}^{2N} n\ a_n\ \left(\int\left[\sum_{j=1}^\infty
  \l_j^{-\a}\left(\toznl (e_j)\right)^2\right]^4
\ \mu(dz)\right)^{1/8} \left( \int \varphi^4\ d\mu\right)^{1/8}\\
& < & \infty
\end{eqnarray*}
since $\a > 0$ by \refeq{5.8} (cf. the proof of \cite[Theorem
7.5]{RZ92}). Furthermore, for all $j\in\N$
\begin{eqnarray*}
&&_-(\b (z),\l_j^{-\a/2} e_j)_+ = \left( (-\Delta +1)_N^{-\delta/2}\beta,
\lambda^{\delta/2+\alpha-\alpha/2}\ e_j\right)_{H_0}\\
&=& - \lambda_j^{-\alpha -\delta/2 +1}\ {}_{H_{-\d}}(z,e_j)_{H_\d}\
\lambda^{\delta/2 + \alpha/2}\\ 
&&-\lambda_j^{-\alpha/2-\delta/2} \sum_{n=1}^{2N} n\ a_n\ \toznl (e_j)\
\lambda^{\delta/2}\\
&=& -\lambda_j\ {}_{H_{-\delta}}(z,\lambda_j^{-\alpha/2}e_j)_{H_\delta} -
\sum_{n=1}^{2N} n\ a_n\ \toznl (\lambda_j^{-\alpha/2}e_j) .
\end{eqnarray*}
It
follows by \cite[Proposition 7.2]{RZ92} that \refeq{2.1} holds. Defining for
$m\in\N$
$$
\delta^m (z):=-\sum_{j=1}^m \lambda_j^{-\a-\d/2+1}\ {}_{H_{-\d}}(z,e_j)_{H_\d}\
\l_j^{\delta /2}e_j\ ,
$$
we check immediately that conditions (i)--(iv) of Theorem \ref{theorem 1}
hold with $c_+ =0=c(\eps_0)$ and $\eps_0=1$, since $\l_j\le 0$ $\forall
j\in\N$. Hence the corresponding Dirichlet operator $(-\Delta-\
{}_{-}(\beta,\nabla\cdot)_{+}\restriction \FCbi )$ is essentially self--adjoint
in $L^2(\calH_-,\nu)$.
\bigskip

For probabilistic consequences concerning uniqueness of the corresponding
martingale problem (resp. the associated infinite--dimensional stochastic
differential equation) we refer to \cite{RZ92} or the detailed discussion in
\cite{AR95} (see, in particular, \cite[Theorem 3.5]{AR95}).
\medskip

Finally, we emphasize that the stronger results on uniqueness in $L^p$ 
presented in the next section apply in the above case also.

\section{Further results and final remarks}

So far we have discussed the uniqueness problem in $L^2$ and in $L^1$.
In both cases we used some a--priori estimates. 
In the case  of $L^1$-uniqueness it was a relatively simple estimate
(Proposition 3), whereas in case of $L^2$-uniqueness a much harder estimate
was used (Proposition 4).
The problem is naturally extended to the $L^p$ setting, via Definition
\ref{definition 3}.
A result on  $L^p$-uniqueness can be proved along the same lines as in the proof of
Theorem 1. The only difference is that one has to have a  uniform estimate of
$\nabla u_n$ in the corresponding $L^{2p}$-space, similar to that obtained in Proposition 4.
Also the idea of obtaining such an estimate is the same as in Proposition 4, with one 
exception: one has to change the test function. Namely, the appropriate test function
is $\n_i u_n(\v \n u_n\v^2+\eps)^{p-2},~\eps>0$, instead of $\n_i u_n$ used in the
proof of Proposition 4 after differentiating the equation. We do not present this estimate 
and its proof here because it contains too many tedious technicalities. We refer to \cite{L1}
where a similar estimate was obtained in a somewhat simpler situation. 
Our present result on  $L^p$-uniqueness, a generalization
of Theorem 1, reads as follows.

\begin{theorem}\label{theorem 3}
Let $\b\in L^p,~\b=\a+\d$, $\v \a\v_0\in L^{2p}$, $\v \d\v_-\in L^2$.
Suppose that there exists a sequence of mappings $(\d^m)_{m\in\Bbb N},\ 
\d^m:\Hscr_-\mapsto\Hscr_-, m\in \Bbb N$ such that all conditions of Theorem
\ref{theorem 1} hold.
Then $(\D+(\b,\n\cdot)\restriction \FCbi)$
has a unique extension which generates a $C_0$-semigroup on 
$L^p$
for all $p\in (1+ {1\over 1+\sqrt{\eps_0}}, 1+{1\over 1-\sqrt{\eps_0}})$.
\end{theorem}

\begin{rem}\label{remark 2}
The uniqueness result in \cite{L1} is now obtained  as a particular case of 
Theorem \ref{theorem 3}:
one has to take $\d=0$. In this case $\eps_0=1$ and the interval of uniqueness becomes
$3/2<p<\infty$. We do not know at the moment whether the assumption on $p$ is just a 
technical restriction or it reflects the essence of the problem. 
\end{rem}
The analysis of the proof of Theorem \ref{theorem 2}  shows that it can also 
be extended to the $L^p$
setting, for $p\in [1,2)$.  The corresponding
 result (i.e., Theorem \ref{theorem 4} below) only complements  Theorem
\ref{theorem 3} for the case 
$p\leq 1+ \frac{1}{1+\sqrt{\eps_0}}$. Otherwise, it is, 
of course, contained in the latter result.

\begin{theorem} \label{theorem 4}
Let $1\leq p<2$, $\b=\a+\d,~~\v \a\v_0\in L^{\frac{2p}{2-p}},~~\v \d\v_-\in 
L^2$.
Suppose that there exists a sequence of mappings $(\d^m)_{m\in\Bbb N},\ 
\d^m:\Hscr_-\mapsto\Hscr_-, m\in \Bbb N$ such that

(i) $(\d^m)_{m\in\Bbb N}$ satisfies conditions (i), (iii) of Theorem
\ref{theorem 2}. 

(ii) $\v\d-\d^m\v_-\longrightarrow 0$  in $L^p(\Hscr_-,d\nu)$ as
$m\rightarrow\infty$,

\noindent
Then the operator $(\D+(\b,\n\cdot)_0\restriction
\FCbi )$ has  a unique extension which generates a $C_0$-semigroup
on $L^p$.
\end{theorem}
The proof is an obvious modification of the proof of Theorem \ref{theorem 2}.

\bigskip

Applications of these results have already been discussed in the previous
section. 
\medskip

Finally, we mention that
the method used in this paper can be extended to the case of Dirichlet
forms with variable non-smooth coefficients. It is also important to investigate the
$L^p$-uniqueness problem for the Dirichlet operator perturbed by a singular potential,
especially with a form-bounded negative part (see \cite{LSe2} for abstract results in this
direction). We intend to return to these problems in the future.

\end{document}